\documentclass[aps,pre]{revtex4}
\usepackage{amssymb}
\usepackage{amsmath}
\usepackage{esint}
\usepackage{color,xcolor}
\usepackage{graphicx}
\usepackage{hyperref}

\newcommand{\myiiint}{\mathop{\int\!\!\!\!\!\int\!\!\!\!\!\int}}

\begin{document}

\title{The moment of inertia tensor of an oloid}

\author{Sander G. Huisman}
\affiliation{Physics of Fluids Department, Max Planck Center for Complex Fluid Dynamics, and J. M. Burgers
Centre for Fluid Dynamics, University of Twente, P.O. Box 217, 7500AE Enschede, The Netherlands}

\date{\today}

\begin{abstract}
The oloid is defined as the convex hull of two unit circles in perpendicular planes, each passing through the center of the other. In this paper we derive an analytical expression for the moment of inertia tensor of an oloid with uniform density and confirm the result numerically.
\end{abstract}
\maketitle

\section{Introduction}
\begin{figure}[htpb!]
	\includegraphics[width=0.5\textwidth]{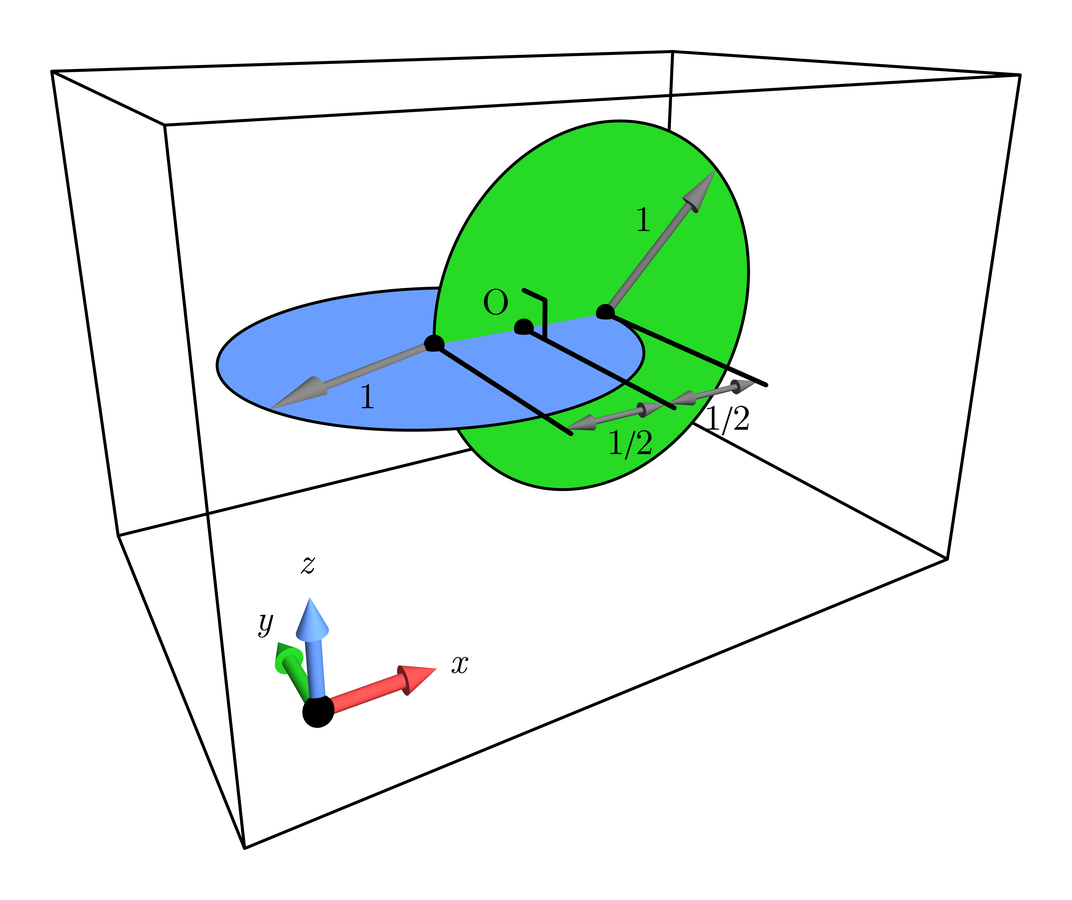}\includegraphics[width=0.5\textwidth]{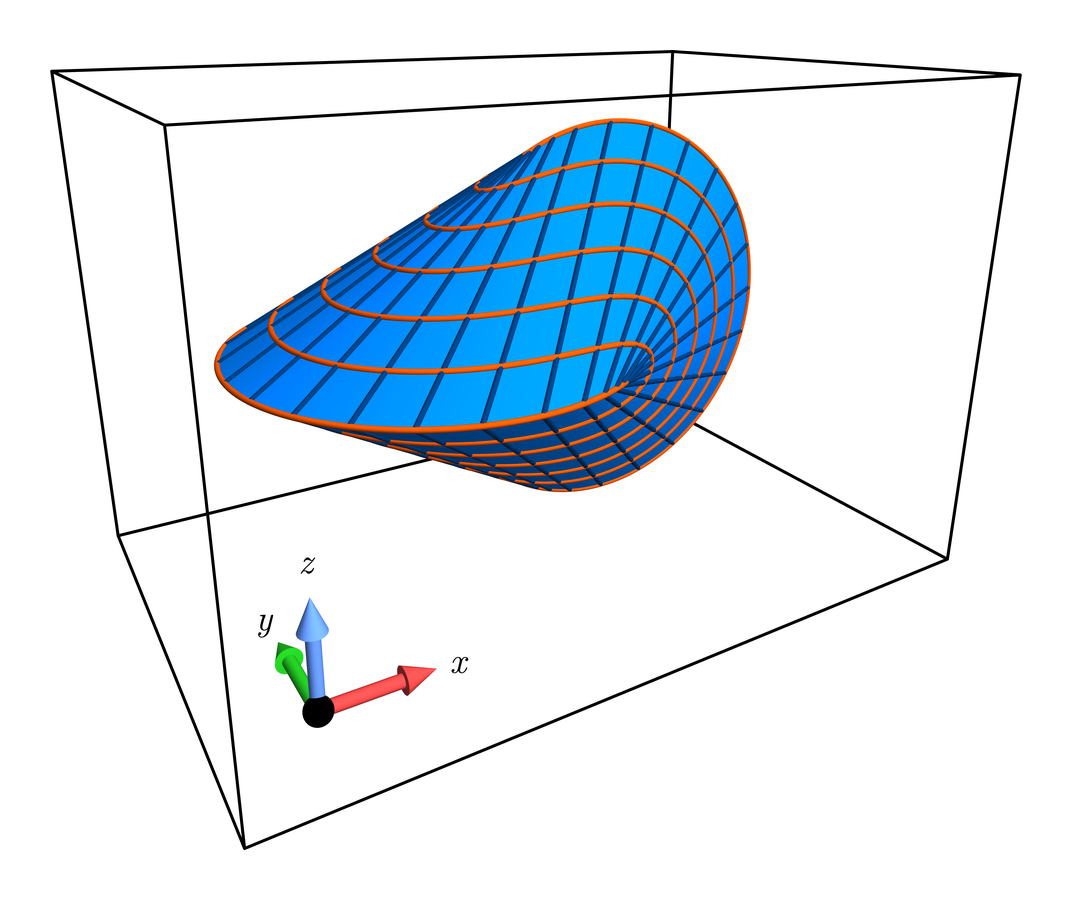}
	\caption{Left: Geometry of the two intersecting circles, here shown as colored disk. The two unit circles lie in perpendicular planes and their centers are located at a distance $1/2$ from the origin O. Right: Convex hull of the configuration on the left, with an added mesh to visualize the 3D surface. The orange mesh-lines are for constant $m$ and the blue mesh-lines are for constant $t$. Both figures include the directions of the axes $x$, $y$, and $z$, see the gnomons in the bottom left, the origin (O) is at the geometrical center of the shape which coincides with the center of mass.}
	\label{fig:OloidGeom}
\end{figure}

The oloid was introduced by Paul Schatz in 1929. It is obtained as the convex hull of two unit circles in perpendicular planes, where each circle passes through the center of the other. The underlying construction and the resulting geometry are shown in Fig.~\ref{fig:OloidGeom}. The oloid is a developable roller, specifically, it is part of the family of so-called two-circle rollers, and it is also closely related to the sphericon, polycons, and platonicons. These types of particles have even been customized to follow a specific path \cite{sobolev2023solid}. The oloid has also been used in the field of fluid mechanics to study the settling of such a particle inside a fluid \cite{flapper2025settling}. 

We place the origin (O) at the midpoint of the two centers, such that the centers are a distance $1/2$ from the origin. The convex hull can then be described parametrically as \cite{dirnbock1997development}:
\begin{align}
	P(m,t) = \begin{pmatrix}
		x(m,t) \\
		y(m,t) \\
		z(m,t) \\
	\end{pmatrix} 
	&=
\begin{pmatrix}
	-\frac{1}{2} + \frac{m}{\cos (t)+1}+(m-1) \cos (t) \\
(1-m) \sin (t) \\
\pm \frac{m \sqrt{2 \cos (t)+1}}{\cos (t)+1}
\end{pmatrix}, \label{eq:para}
\end{align}
where $0 \leq m \leq 1$ and $-2\pi/3 \leq t \leq 2\pi/3$. Here the $\pm$ describes the upper and lower part of the hull; combining the two then gives the complete hull. Previously, it has been shown that the area of the object is $A=4\pi$ \cite{dirnbock1997development}. That paper also reported the numerical value of the volume. Here we express the result as the sum of two elliptic integrals (which can be found by symbolically integrating their Eq.~32):
\begin{align}
 V &=  \frac 23 K\left(\frac{3}{4}\right) +\frac 43 E\left(\frac{3}{4}\right) = 3.05241846842437\ldots
\end{align}
where $K$ and $E$ are the complete elliptic integral of the first kind and second kind, respectively:
\begin{align}
	K(m) &= \int \limits_{0}^{\pi/2} \frac{1}{\sqrt{1-m \sin(\theta)^2}} \mathrm{d}\theta \\ 
	E(m) &= \int \limits_{0}^{\pi/2} \sqrt{1-m \sin(\theta)^2} \mathrm{d}\theta.
\end{align}
Modeling the rotational dynamics of an oloid requires the moment of inertia tensor, which is composed of 9 components:
\begin{align}
 I_{ij} &= \myiiint \limits_V \rho \left(x_k x_k \delta_{ij} - x_i x_j \right) \mathrm{d}V = \begin{pmatrix} I_{xx} & I_{xy} & I_{xz} \\ I_{yx} & I_{yy} & I_{yz} \\ I_{zx} & I_{zy} & I_{zz}  \end{pmatrix} 
\end{align}
where we used the Einstein summation convention, and we will assume $\rho=1$ from now on (uniform unit density). We know that the tensor is symmetric, such that we are looking for (at most) 6 unique elements (3 diagonal, 3 off-diagonal). We choose our axes as in figure \ref{fig:OloidGeom}, the origin (O) is aligned with the center of mass. The axis defined by the intersection of the two symmetry planes of the oloid (the $x$ direction through the origin O) is a principal axis. The remaining two axes can be chosen freely (provided they are mutually perpendicular and perpendicular to the first axis). For convenience we align them with the planes of the circles. By choosing the principal axes the moment of inertia tensor is simplified and all the off-diagonal terms are zero:
\begin{align}
 I_{ij} &= \begin{pmatrix} I_{xx} & 0 & 0 \\ 0 & I_{yy} & 0 \\ 0 & 0 & I_{zz}  \end{pmatrix}.
\end{align} 
Alternatively, one can see from the definition of the moment of inertia (for uniform unit density):
\begin{align}
	I_{xy} &= \myiiint \limits_V x y \mathrm{d}V
\end{align}
that because the body is symmetric with respect to $y\rightarrow -y$, this transformation implies:
\begin{align}
	I_{xy} &= -I_{xy}
\end{align}
from which we conclude that the value must be 0. We can do this also for the $I_{yz}$ term, as well as for the $I_{xz}$ term (because of the symmetry in the $z$ plane), and because the tensor is symmetric thus also for the $I_{yx}$, $I_{zy}$, and $I_{zx}$ terms.

In addition, because of the additional symmetry of the particle we see that $I_{yy}=I_{zz}$ such that we are left to calculate two values:
\begin{align}
	I_{xx} &= \myiiint \limits_V \left(y^2+z^2\right) \mathrm{d}V\\
	I_{yy} &= \myiiint \limits_V \left(x^2+z^2\right) \mathrm{d}V.
\end{align}
To evaluate these integrals we use the divergence theorem to convert the volume integral to a surface integral:
\begin{align}
	I_{xx} &= \oiint \limits_{S} \begin{pmatrix} x\left(y^2+z^2\right)\\0\\0 \end{pmatrix}\cdot \mathrm{d}\vec S \\ 
	I_{yy} &= \oiint \limits_{S} \begin{pmatrix} 0 \\ y\left(x^2+z^2\right)\\0 \end{pmatrix}\cdot \mathrm{d}\vec S
\end{align}
where $S$ is the surface of $V$ and $\mathrm{d}\vec S$ is the infinitesimal surface area multiplied by the local outward-facing normal vector. Using the parametric representation of the surface, this becomes:
\begin{align}
	I_{xx} &= \oiint \limits_{S} \begin{pmatrix} x\left(y^2+z^2\right)\\0\\0 \end{pmatrix}\cdot \mathrm{d}\vec S = 2 \int \limits_0^1 \int \limits_{-2\pi/3}^{2\pi/3} \begin{pmatrix} x\left(y^2+z^2\right)\\0\\0 \end{pmatrix}\cdot \left(\frac{\mathrm{d}P}{\mathrm{d}m} \times \frac{\mathrm{d}P}{\mathrm{d}t} \right) \mathrm{d}t\mathrm{d}m \label{eq:intixx}\\
	I_{yy} &= \oiint \limits_{S} \begin{pmatrix} 0 \\ y\left(x^2+z^2\right)\\0 \end{pmatrix}\cdot \mathrm{d}\vec S = 2\int \limits_0^1 \int \limits_{-2\pi/3}^{2\pi/3} \begin{pmatrix} 0 \\ y\left(x^2+z^2\right)\\0 \end{pmatrix}\cdot \left(\frac{\mathrm{d}P}{\mathrm{d}m} \times \frac{\mathrm{d}P}{\mathrm{d}t} \right) \mathrm{d}t\mathrm{d}m \label{eq:intiyy}
\end{align}
where the factor of 2 accounts for the upper and lower parts of the surface. We now calculate the last bracketed term; the cross product of the surface tangents (here we have take parameterisation of the top shell; the $+$ version of the $z$ component of Eq.~\ref{eq:para}):
\renewcommand*{\arraystretch}{1.5}
\begin{align}
	\frac{\mathrm{d}P}{\mathrm{d}m} &= \begin{pmatrix}  \cos (t)+\frac{1}{\cos (t)+1} \\
 -\sin (t) \\
 \frac{\sqrt{2 \cos (t)+1}}{\cos (t)+1}  \end{pmatrix} \\
	\frac{\mathrm{d}P}{\mathrm{d}t} &= \begin{pmatrix}
 \sin (t) \left(m \left(\frac{1}{(\cos (t)+1)^2}-1\right)+1\right) \\
 (1-m) \cos (t) \\
 \frac{m \sin (t) \cos (t)}{(\cos (t)+1)^2 \sqrt{2 \cos (t)+1}} 
\end{pmatrix}.
\end{align}
We then have:
\begin{align}
 \frac{\mathrm{d}P}{\mathrm{d}m} \times \frac{\mathrm{d}P}{\mathrm{d}t} &= \begin{pmatrix}
 \frac{\cos (t) ((3 m-2) \cos (t)-1)}{2\cos\left(\frac{t}{2}\right)^2 \sqrt{2 \cos (t)+1}} \\
 \frac{\sin (t) ((2-3 m) \cos (t)+1)}{(\cos (t)+1) \sqrt{2 \cos (t)+1}} \\
 \frac{(2-3 m) \cos (t)+1}{\cos (t)+1} \label{eq:cross}
\end{pmatrix}.
\end{align}
We now have all the ingredients to evaluate the integrals of Eqs.~\ref{eq:intixx} and \ref{eq:intiyy}. Note that $x$, $y$, and $z$ are defined in terms of $m$ and $t$, see Eq.~\ref{eq:para}. After substituting everything we obtain:
{\small{
\begin{align}
	I_{xx} &= \int \limits_0^1 \int \limits_{-2\pi/3}^{2\pi/3} \frac{\cos (t)  ((3 m-2) \cos (t)-1) \left(\frac{m}{\cos (t)+1}+(m-1) \cos (t)-\frac{1}{2}\right) \left(\frac{m^2 (2 \cos (t)+1)}{(\cos (t)+1)^2}+(m-1)^2 \sin(t)^2\right)}{\cos\left(\frac{t}{2}\right)^2\sqrt{2 \cos (t)+1}}  \mathrm{d}t\mathrm{d}m \label{eq:finaleqxx}\\
	I_{yy} &= \int \limits_0^1 \int \limits_{-2\pi/3}^{2\pi/3} \frac{2 (1-m) \sin(t)^2 ((2-3 m) \cos (t)+1) \left((m-1)^2 \cos(t)^2-(m-1) \cos (t)+m \left(2 m+\frac{1}{\cos (t)+1}-2\right)+\frac{1}{4}\right)}{(\cos (t)+1) \sqrt{2 \cos (t)+1}} \mathrm{d}t\mathrm{d}m. \label{eq:finaleqyy}
\end{align}
}}
The integration with respect to $m$ is polynomial of degree 4 and results in:
{\small{
\begin{align}
	I_{xx} &= \int \limits_{-2\pi/3}^{2\pi/3} -\frac{\cos(t)^2 (-2510 \cos (t)+547 \cos (2 t)+1129 \cos (3 t)+648 \cos (4 t)+181 \cos (5 t)+21 \cos (6 t)-1744) }{15360 \cos\left(\frac{t}{2}\right)^8\sqrt{2 \cos (t)+1}} \mathrm{d}t \\
	I_{yy} &= \int \limits_{-2\pi/3}^{2\pi/3} \frac{(542 \cos (t)+322 \cos (2 t)+122 \cos (3 t)+21 \cos (4 t)+361) \tan\left(\frac{t}{2}\right)^2}{240 \sqrt{2 \cos (t)+1}} \mathrm{d}t.
\end{align}
}}
These expressions can be integrated with respect to $t$ manually or using a computer algebra system (CAS):
\begin{align}
	I_{xx} &= \frac{32}{45} E\left(\frac{3}{4}\right)-\frac{2}{45} K\left(\frac{3}{4}\right) = 0.76535025749314262939\ldots \\
	I_{yy} = I_{zz} &= \frac{71}{45} E\left(\frac{3}{4}\right)-\frac{19}{90} K\left(\frac{3}{4}\right) = 1.45551287346920034498\ldots
\end{align}
To calculate these constants precisely one can use the link between elliptic integrals and Gauss's arithmetic-geometric mean (AGM) to develop an algorithm that converges quadratically. We have confirmed the value of these numbers by numerically integrating Eqs.~\ref{eq:finaleqxx} and \ref{eq:finaleqyy}. In addition, a Monte Carlo simulation in which $10^6$ random points are generated inside the hull and calculating the moment of inertia also confirms these numbers up to 2 decimal places. For an oloid constructed from circles with a radius $r$ and unit density the moment of inertia is given by:
\begin{align}
 I(r) &= \begin{pmatrix} I_{xx} & 0 & 0 \\ 0 & I_{yy} & 0 \\ 0 & 0 & I_{zz}  \end{pmatrix} r^5.
\end{align} 
\begin{acknowledgments}
We thank Mees Flapper, Detlef Lohse, Bernhard Mehlig, John Sader, Federico Toschi, Leen van Wijngaarden, Xander de Wit, and Greg Voth for stimulating discussions about the oloid.
\end{acknowledgments}

\section*{Appendix: Area calculation}
Finding the area of the oloid can be found by integrating the magnitude of the normal vector from Eq.~\ref{eq:cross}:
\begin{align}
A &= 2\int \limits_0^1 \int \limits_{-2\pi/3}^{2\pi/3}  \left|\left| \frac{\mathrm{d}P}{\mathrm{d}m} \times \frac{\mathrm{d}P}{\mathrm{d}t} \right|\right| \mathrm{d}t\mathrm{d}m \\
A &= \int \limits_0^1 \int \limits_{-2\pi/3}^{2\pi/3} \frac{(4-6 m) \cos (t)+1}{\cos\left(\frac{t}{2}\right)\sqrt{2 \cos (t)+1}} \mathrm{d}t\mathrm{d}m \\
A &= \int \limits_{-2\pi/3}^{2\pi/3} \frac{\cos (t)+2}{\cos\left(\frac{t}{2}\right) \sqrt{2 \cos (t)+1}} \mathrm{d}t\\
A &= 4\pi
\end{align}
confirming the results of Ref.~\cite{dirnbock1997development}.
\end{document}